\documentclass[12pt]{amsart}
\usepackage{amssymb}
\usepackage{fullpage}
\usepackage[colorlinks=true]{hyperref}

\newcommand{\bbR}{\mathbb{R}}

\newtheorem*{theorem-nonum}{Theorem}

\newtheorem*{conjecture-nonum}{Conjecture}

\theoremstyle{definition}

\newtheorem*{remark-nonum}{Remark}
\newtheorem*{remarks-nonum}{Remarks}

\begin{document}

\title{On the Guy-Kelly Conjecture for the No-Three-In-Line Problem}

\author{Paul M Voutier}
\address{London, UK}
\email{Paul.Voutier@gmail.com}



\begin{abstract}
We provide details of the error Gabor Ellmann found in 2004 in the heuristic argument of
Guy and Kelly in \cite{G-K}. This led to a correction of their conjectured upper
bound for the no-three-in-line problem. However, details of the issue and its
correction, including the actual location of the issue, while simple, do not
seem to have appeared in the literature previously. That said, very recent
work of Prellberg \cite{Pr} does contain a derivation of the corrected
conjectured upper bound.
\end{abstract}

\maketitle

\section{Introduction}

Following the notation in \cite{G-K}, we let $S_{n}$ be the set of $n^{2}$ points
in $\bbR^{2}$ with integer coordinates $(x, y)$,
$1 \leq x, y < n$. Let $f_{n}$ be the maximum size of a subset $T$ of
$S_{n}$ such that no three points of $T$ are collinear. The study of $f_{n}$ is
the well-known \emph{no-three-in-line problem} in discrete geometry.

The upper bound $f_{n} \leq 2n$
is an easy application of the pigeon-hole principle.
In fact, it is known that for $n \leq 64$, we have $f_{n}=2n$ \cite{Flam1}. The
extension beyond $n=60$ to $n=64$ having only been attained in mid-February 2026
by Prellberg.

However, it is conjectured that for $n$ sufficiently large, $f_{n}<2n$. In their
paper \cite{G-K}, Guy and Kelly presented a heuristic that led them to conjecture that
\begin{equation}
\label{eq:gk-1}
f_{n} \sim \left( 2\pi^{2}/3 \right)^{1/3} n.
\end{equation}

In private communication in March 2024, Gabor Ellmann pointed out an error to
R.K. Guy (see the Comments section of \cite{OEIS}, as well as \cite{MP}):
\begin{verbatim}
From R. K. Guy, Oct 22 2004: (Start)
...
As recently as last March, Gabor Ellmann pointed out an error in our heuristic
reasoning, which, when corrected, gives 3*c^2 = Pi^2, or c ~ 1.813799. (End)
\end{verbatim}

As stated there by Guy, Ellman's correction leds to the following conjecture.

\begin{conjecture-nonum}
\begin{equation}
\label{eq:pv-1}
f_{n} \sim \left( \pi/\sqrt{3} \right) n.
\end{equation}
\end{conjecture-nonum}

Guy stated his intention in sending a correction to Can. Math. Bull in private
communication to Ed Pegg in Oct 2004 \cite{MP} but never did so before his death
in 2020.

To the best of our knowledge, the precise nature of this error, its correction
and how it leads to \eqref{eq:pv-1} has never been published.
The nearest explanation found online is the following:
\begin{verbatim}
In March 2005 Gabor Ellmann discovered an error in Guy & Kellys published
derivation of 1968. 
If this geometric/algebraic transformation-mistake is corrected it results:
asymptotically are only PI/3^1/2 n = 1.814... n points selectable. 
\end{verbatim}
This appears on A. Flammenkamp's website, \cite{Flam2}. Note that he gives the
date as March 2005. We have chosen to follow the date from \cite{MP} and \cite{OEIS}.

Our sole contribution here is provide this correction in the literature. Although
the error is very localised, as described in our final remark below, for completeness,
as well as to show the correction, we have included the entire Guy-Kelly argument below.

As noted in our abstract, very recently, Prellberg gave a derivation of the
corrected version of the conjecture. See pages~2--3 of \cite{Pr}.

Before moving on to Guy and Kelly's heuristic argument, we note that recently
some people have questioned the independence assumption in their argument (we use
it below to obtain \eqref{eq:independence}). See, for instance,
\url{https://11011110.github.io/blog/2018/11/10/random-no-three.html} discussing
a seminar of Nathan Kaplan on this topic in November 2018. Also some empirical
evidence such as from $90$-degree rotationally symmetric grids (for $n$ even),
calls this into question. Prellberg's work cited above (see the last few paragraphs
of Section~1 on page~3 of \cite{Pr}) also contains a discussion of this and
doubts about the validity of their heuristic argument. We thank Thomas Prellberg
for these references and observations.

\section{The Heuristic Argument}

In their paper \cite{G-K}, Guy and Kelly prove the following.

\begin{theorem-nonum}
The number, $t_{n}$, of sets of $3$ collinear points that can be chosen from
$S_{n}$ is
\[
t_{n}=\frac{3}{\pi^{2}} n^{4} \log(n) + O \left( n^{4} \right).
\]
\end{theorem-nonum}

After proving this result, they procced use it to give a probabilistic argument for
\eqref{eq:gk-1}.

The probability that three points, chosen at random, should be in line is thus
\[
\left( \frac{3}{\pi^{2}} n^{4} \log(n) + O \left( n^{4} \right) \right) / \binom{n^{2}}{3}
= \frac{18\log(n)}{\pi^{2}n^{2}} + O \left( \frac{1}{n^{2}} \right)
\]
and so the probability that three such points should not be in line is
\[
1-\frac{18\log(n)}{\pi^{2}n^{2}} + O \left( \frac{1}{n^{2}} \right).
\]

If we assume that the events are independent, the probability that, for any
fixed $k>0$, a set of $kn$ points contains no three in line is
\begin{align}
\label{eq:independence}
\left( 1-\frac{18\log(n)}{\pi^{2}n^{2}} + O \left( \frac{1}{n^{2}} \right) \right)
^{\binom{kn}{3}}
&=\exp \left( -\frac{3k^{3}n \log(n)}{\pi^{2}}+ O(n) \right).
\end{align}
This equality holds because for positive integers $m$, we have
$(1+x)^{m}=\exp \left( mx +O \left( mx^{2} \right) \right)$,
provided that $|x|$ is small. Here we have $m={\binom{kn}{3}}
=k^{3}n^{3}/6+O \left( n^{2} \right)$ and $x=-18\log(n)/ \left( \pi^{2}n^{2} \right)
+O \left( n^{-2} \right)$, so $mx=-3k^{3}n \log(n)/\pi^{2}+ O(n)$ and
$mx^{2}=O \left( \log^{2}(n)/n \right)$.

Hence an estimate of the number of solutions to the no-three-in-line problem with
$kn$ points is given by
\[
\binom{n^{2}}{kn} n^{-3k^{3}n/\pi^{2}}
\cdot e^{O(n)}.
\]

As $n \rightarrow \infty$, we have $\left( n^{2} \right)!/ \left( n^{2}-kn \right)!
\sim \left( n^{2} \right)^{kn}$. We also apply Stirling's approximation to $(kn)!$ and hence obtain
\begin{equation}
\label{eq:5}
O \left( n^{kn-3k^{3}n/\pi^{2}}\cdot e^{O(n)} \right),
\end{equation}
as the estimate for the number of solutions. This is the analogue of equation~(5)
in \cite{G-K}.

If $k-3k^{3}/\pi^{2}<0$, then this estimate goes to $0$ as $n \rightarrow
\infty$. This holds provided that $k> \pi/3^{1/2}$, completing the heuristic
argument for the conjecture.

\begin{remark-nonum}
We see that the precise location of the error in \cite{G-K} is on line~-4 of page~530,
where $-2+3k^{3}/\pi^{2}$ there should read $-k+3k^{3}/\pi^{2}$.

This appears to have arisen by using $2n$ instead of $kn$ in the estimates they
used to obtain their equation~(5). This is akin to estimating
$\left( n^{2} \right)!/ \left( n^{2}-2n \right)!$ and $(2n)!$ instead of
$\left( n^{2} \right)!/ \left( n^{2}-kn \right)!$ and $(kn)!$, respectively, in
our argument for our equation~\eqref{eq:5} above.
\end{remark-nonum}

\end{document}